\documentclass{amsart}
\usepackage{amsthm,amssymb,amsmath}

\theoremstyle{definition}

\newtheorem{theorem}{Theorem}[section]
\newtheorem{lemma}[theorem]{Lemma}

\newtheorem{example}[theorem]{Example}

\theoremstyle{remark}

\numberwithin{equation}{section}
\newcommand{\supp}{\mbox{supp}}

\newcommand{\const}{\mbox{const}}

\newcommand{\cl}{\mbox{cl}}

\newcommand{\loc}{{\rm loc}}

\newcommand{\mes}{\mbox{mes}}

\begin{document}

\title[The space of distributions with discontinuous test functions]
{The space of distributions with discontinuous test functions and a family of zero-sum games with discontinuous payoffs}

\author[V.Derr] {V.~Derr}

\address{Faculty of Mathematics, Udmurtia State University, Universitetskaya St., 1 (Bld. 4), Izhevsk, 426034, Russia}

\email{derr@uni.udm.ru}

\author[D.Kinzebulatov] {D.~Kinzebulatov}

\address{Department of Mathematics, University of Toronto, Toronto, Ontario, Canada M5S 2E4}

\email{dkinz@math.toronto.edu}



\begin{abstract}
In the present paper we consider one class of zero-sum games with discontinuous payoffs which may have no solutions in the sets of pure or mixed strategies. We show that, however, the solution always exists in the set of so-called $\mathcal R'$-mixed strategies which are the elements of the space $\mathcal R'$ of distributions with discontinuous test functions.
The advantages of the new space of distributions (in comparison with the classical space $\mathcal D'$ of distributions with continuous or smooth test functions), that is, the possibility to define in $\mathcal R'$ the operations of integrations of distributions and multiplication of distributions by discontinuous functions, which are correct in the sense that they are everywhere defined, continuous and coincide with the ordinary operations for regular distributions, are crucial for our proof of existence of solution.
\end{abstract}

\maketitle
%

\section{Introduction}
For the past decades the progress of L.~Schwartz' distribution theory was highly motivated by efforts to overcome the following well-known insufficiencies of the classical space $\mathcal D'(\Omega)$ of distributions with the test functions continuous  and having compact support in $\Omega \subset \mathbb R^n$, $\Omega$ is an open set: the impossibility to define in $\mathcal D'(\Omega)$ the correct operation of multiplication of distributions by discontinuous functions (or,  more generally, of multiplication of two distributions) as well as the correct operation of integration of distributions, where, following \cite{Shi}, we say that an operation is \textit{correct} if it is defined for all distributions, continuous and coincides with the ordinary one for regular distributions (see \cite{Bag1,Bag2,Sar,Sar2,Sar3,Ses} and further references therein).
The numerous applications of distribution theory to ordinary and partial differential equations (e.g., see \cite{Bag2,ColM3,Fil,Sar2,Sar3}), where the necessity to integrate distributions and to multiply distributions by discontinuous functions arise, demonstrate the importance of these operations. 
In the present paper we define certain extension of the classical space of distributions $\mathcal D'(\Omega)$ where the correct operations of integration and multiplication by discontinuous functions are defined, and use this extension to provide the existence of Nash equilibrium for a class of zero-sum games with discontinuous payoffs.

As it follows from the classical definition of the product of a distribution $f$ and an ordinary function $g$ as the distribution given by
\begin{equation}
\label{inprod}
(gf,\varphi) = (f,g\varphi),
\end{equation}
where $\varphi$ stands for a test function,
the ordinary function $g$ and the test function $\varphi$ must possess the same degree  of regularity in order for the product $gf$ to be defined. In particular, in the space $\mathcal D'(\Omega)$ the ordinary multiple $g$ must be continuous. 
This observation suggests the idea of extension of the classical space of continuous test functions $\mathcal D(\Omega)$ so that the test functions are allowed to be discontinuous, where here and below under discontinuous functions we understand the elements of the closure of the algebra of piecewise-constant functions in the topology of uniform convergence, called regulated functions (see \cite{Dieu,Dav}). 
This idea appeared first in \cite{Kur,Kur2} where the space of distributions with the test functions defined on $\mathbb R$ and possibly having first-kind discontinuity at the origin was introduced.
In what follows, we consider the space $\mathcal R'(\Omega)$ of distributions defined on the space $\mathcal R(\Omega)$ of discontinuous test functions having compact support in $\Omega$ and possessing, in a sense, an arbitrary set of points of discontinuity. We employ equality (\ref{inprod}) to define a continuous, commutative and associative operation of multiplication of distributions by discontinuous functions, and define continuous operation of integration of distributions by the formula
\begin{equation*}
\int_S f dx = (f,\chi_S), 
\end{equation*}
where $\chi_S \in \mathcal R(\Omega)$ stands for the characteristic function of a bounded subset $S \Subset \Omega$.
In what follows, we show that any distribution in $\mathcal D'(\Omega)$ admits a continuous linear extension to the space of discontinuous test functions $\mathcal R(\Omega)$. Let us note that since all elements of the space $\mathcal D'(\Omega)$ are measure-type distributions, the elements of our distribution space $\mathcal R'(\Omega)$ are also, in a sense, of measure-type.
We further describe the family of delta-functions in $\mathcal R'(\Omega)$ concentrated at a given point: for $n=1$ the classical delta-function $\delta \in \mathcal D'(\mathbb R)$ extends to a family of delta-functions $\delta^\alpha \in \mathcal R'(\mathbb R)$, where $\alpha$ is a function,
\begin{equation*}
\alpha:\{-1,1\} \to \mathbb R, \quad \alpha(1)+\alpha(-1)=1,
\end{equation*}
so that $\delta^\alpha$ becomes an affine combination of the \textit{right} and the \textit{left} delta-functions, respectively,
\begin{equation*}
\delta^\alpha=\beta\delta^{+} +(1-\beta)\delta^{-}, \quad \beta=\alpha(1) \in \mathbb R, \quad (\delta^+,\varphi)=\varphi(0+), \quad (\delta^-,\varphi)=\varphi(0-),
\end{equation*}
$\varphi$ is the test function,
and the product of $\delta^\alpha$ with Heaviside function $\theta=\chi_{(\infty,0)}$ is given by
\begin{equation}
\label{intro_delta}
\theta\delta^\alpha=\beta\delta^+.
\end{equation}
Due to the continuity of the operation of multiplication in $\mathcal R'$ the same result can be obtained if the delta-function in (\ref{intro_delta}) is replaced by terms of the corresponding delta-sequence.

The main purpose of our extension is to provide the existence of solution (i.e., Nash equilibrium) for a family of zero-sum games with discontinuous payoff functions, which in general do not have a solution in the set of pure strategies or in the set of mixed strategies, yet possess a solution in the set of so-called $\mathcal R'$-mixed strategies, which are the elements of the space $\mathcal R'(\Omega)$. For example,
the following zero-sum game $G=(X_1,X_2,\rho)$ with the sets of pure strategies $X_1=X_2=(-1,1)$ and the payoff function given by
\begin{equation*}
\rho(x_1,x_2)=\left\{
\begin{array}{ll}
1,&x_1,x_2>0,x_1+x_2<1 \text{ or } \\
&x_1,x_2<0,x_1+x_2>-1, \\
0, & \text{otherwise}
\end{array}
\right.
\end{equation*}
does not have solution in the sets of pure strategies or classical mixed strategies. However, as it turns out, there exists solution of $G$ in the set of $\mathcal R'$-mixed strategies, which is the pair
\begin{equation*}
0.5(\delta^++\delta^-), \quad 0.5(\delta^++\delta^-),
\end{equation*}
as it follows from a general result (Theorem \ref{gameteo}) proved below. The continuity of the operations of multiplication and integration allows us to provide certain probabilistic interpretation for $\mathcal R'$-mixed strategies.

\section{Regulated functions}

\textbf{1.~}Let $\Omega \subset \mathbb R$ be an open interval.
Following \cite{Dieu}, we call function $g:\Omega \mapsto \mathbb R$ regulated, if it is bounded and for every $x \in \cl(\Omega)$ there exist one-sided limits
$g(x+)$, $g(x-)$.
We denote the algebra of regulated functions by $\hat{\mathbb G}(\Omega)$ and endow with the standard sup-norm.

The function $g \in \hat{\mathbb G}(\Omega)$ is called
piecewise-constant, if for any bounded subset $\Gamma \subset \Omega$ the image $g(\Gamma)$ is finite.
In what follows, we denote the algebra of piecewise-constant functions by $\hat{\mathbb{PC}}(\Omega)$.

\begin{lemma}[\cite{Dieu}]
\label{lem0}
$\hat{\mathbb{PC}}(\Omega)$ is dense in $\hat{\mathbb G}(\Omega)$.
\end{lemma}

Let us define $$J=\{g \in \hat{\mathbb G}(\Omega):g(x+)=g(x-)=0~(x \in \Omega)\};$$ $J$ is a closed ideal in $\hat{\mathbb G}(\Omega)$, so we may define
$\mathbb G(\Omega)=\hat{\mathbb G}(\Omega)/J$.
The norm in $\mathbb G(\Omega)$ is given by
\begin{equation}
\|g\|_{\mathbb G(\Omega)}=\sup_{x \in \Omega}\max\{|g(x+)|,|g(x-)|\}.
\end{equation}

\begin{lemma}
\label{banach}
$\mathbb G(\Omega)$ is a Banach algebra.
\end{lemma}
\begin{proof}
$\hat{\mathbb G}(\Omega)$ is complete according to \cite{Dieu}. The ideal $J$ is closed in $\hat{\mathbb G}(\Omega)$, hence the factor-algebra $\mathbb G(\Omega)$ is also complete \cite{Gam} .
\end{proof}

Given $g \in \mathbb G(\Omega)$, we define $$T(g)= \{x \in \Omega: g(x+) \ne g(x-)\}$$ to be the set of point of discontinuity of $g$.
Similarly, we denote by $\mathbb{PC}(\Omega)\hookrightarrow \mathbb G(\Omega)$ the image of the algebra of piecewise-constant functions under the canonical map $\hat{\mathbb G}(\Omega) \to \mathbb G(\Omega)$.

Let $\mathbb C(\Omega) \hookrightarrow \mathbb G(\Omega)$ be the algebra of continuous elements of $\mathbb G(\Omega)$. 

\begin{lemma}
\label{teo1}
$\mathbb{PC}(\Omega)$ is dense $\mathbb G(\Omega)$.
\end{lemma}
\begin{proof}
Let $\hat{g} \in \hat{\mathbb G}(\Omega)$ be given. Suppose that $g$ is the corresponding equivalence class of $\hat{g}$ in $\mathbb G(\Omega)$. Then $\|g\|_{\mathbb G(\Omega)} \leqslant \|\hat{g}\|_{\hat{\mathbb G}(\Omega)}$, so it suffices to apply Lemma \ref{lem0} to complete the proof.
\end{proof}

\begin{lemma}
Given $g \in \mathbb G(\Omega)$, the set $T(g)$ is at most countable.
\end{lemma}
\begin{proof}
Let us denote by $\hat{T}(\hat{g})$ the set of points of discontinuity of a representative $\hat{g} \in \hat{\mathbb G}(\Omega)$. Then $T(g) \subset \hat{T}(\hat{g})$. Since $\hat{T}(\hat{g})$ is at most countable \cite{Dieu},
$T(g)$ is also at most countable.
\end{proof}

\textbf{2.~}We proceed to the definition of the algebra of regulated functions defined on an open set in $\mathbb R^n$.
Let $\Omega \subset \mathbb R^n$ be an open set,
and let $\mathcal F$ be the family of
finite unions and finite differences of convex subsets of $\Omega$. 
The family $\mathcal F$ is called the appropriate family (in terminology of \cite{Dav}). 
Following \cite{Dav}, we call function $g:\Omega \mapsto \mathbb R$ regulated if it is bounded and for every $x \in \cl(\Omega)$ and for any $\varepsilon>0$ there exist a neighbourhood $U_x=U_x(\varepsilon) \in \mathcal F$ and a family $\{S_i\}_{i=1}^{m} \subset \mathcal F$ such that $U_x=\cup_{i=1}^m S_i$, and
\begin{equation*}
|g(y_1)-g(y_2)|<\varepsilon
\end{equation*}
for all $y_1,y_2 \in S_i$
\rm ($1 \leqslant i \leqslant m$). Analogously, we denote the algebra of regulated functions by $\hat{\mathbb G}(\Omega)$ and endow it with sup-norm. 

The regulated function $g \in \hat{\mathbb G}(\Omega)$ is called piecewise-constant, if for any bounded open subset $\Gamma \subset \Omega$ the restriction $g|_{\Gamma}$ is a linear combination of the characteristic functions $\chi_S$, where $S \in \mathcal F \cap \Gamma$. We denote the algebra of piecewise-constant functions by $\hat{\mathbb{PC}}(\Omega)$.

\begin{lemma}[\cite{Dav}]
\label{lem1}
$\hat{\mathbb{PC}}(\Omega)$ is dense in $\hat{\mathbb G}(\Omega)$.
\end{lemma}

Let $g \in \hat{\mathbb G}(\Omega)$ and $x_0 \in \Omega$ be fixed, let $S^{n-1}$ be the unit sphere in $\mathbb R^n$ centered at $0$. 
For each $s \in S^{n-1}$ we define
\begin{equation}
\label{limit}
g(x_0)(s) = \lim_{t \to 0+} g(x_0+ts).
\end{equation}
Let us show that the function $g(x_0):S^{n-1} \mapsto \mathbb R$ is defined everywhere on $S^{n-1}$. Let $s \in S^{n-1}$ be given. We denote $$L=\{x_0+ts:0<t<1\}, \quad \mathcal F_L = \{A \cap L: A \in \mathcal F\}.$$ Now as it follows from the definition of the family $\mathcal F$, $\mathcal F_L$ is an appropriate family for the line segment $L$, so in virtue of Lemma \ref{lem0} the limit (\ref{limit}) exists.

In what follows, we call $g(x):S^{n-1} \mapsto \mathbb R$ the surrounding value of $g$ at $x$.
Note that if $g$ is continuous at $x$, then by definition $g(x)(\cdot) \equiv g(x)$.

\begin{lemma}
\label{lem5}
For any $x \in \Omega$ the surrounding value $g(x)(\cdot)\in \mathbb L^{\infty}(S^{n-1})$.
\end{lemma}
\begin{proof}
Suppose that $x \in \Omega$ is given. Clearly, $g(x)(\cdot)$ is bounded. Suppose that $\{g_k\}_{k=1}^\infty \subset \hat{\mathbb{PC}}(\Omega)$, $g_k \to g$ uniformly on $\Omega$. 
Then $g_k(x)(\cdot)$ is Lebesgue measurable on $S^{n-1}$. Due to uniform convergence 
we may change the order of limits, so
\begin{equation}
\label{lem6conv}
\sup_{s \in S^{n-1}}(g_k(x)(s)-g(x)(s)) \to 0. 
\end{equation}
Consequently, $g(x)(\cdot)$ is Lebesgue measurable on $S^{n-1}$ as a limit of the uniform convergent sequence of Lebesgue measurable functions. As a result, $g(x)(\cdot)\in \mathbb L^{\infty}(S^{n-1})$.
\end{proof}

We define $$J=\{g \in \hat{\mathbb G}(\Omega): g(x)(\cdot)=0~(x \in \Omega)\}$$ -- a closed ideal in $\hat{\mathbb G}(\Omega)$, and put
$\mathbb{G}(\Omega)=\hat{\mathbb G}(\Omega)/J$.
The elements of the factor-algebra $\mathbb G(\Omega)$ do not have values at points in $\Omega$, but instead possess surrounding values introduced above.

The algebra $\mathbb G(\Omega)$ is endowed with the norm
\begin{equation}
\|g\|_{\mathbb G(\Omega)}=\sup_{x \in \Omega}\{\|g(x)(\cdot)\|_{\mathbb L^{\infty}(S^{n-1})}\}.
\end{equation}

Clearly, in the case when $n=1$ our definitions coincides with the previous ones.

\begin{lemma}
\label{Gbanach}
$\mathbb G(\Omega)$ is a Banach algebra.
\end{lemma}
\begin{proof}
Since $\mathbb G(\Omega)=\hat{\mathbb G}(\Omega)/J$, where $J$ is a closed ideal, and
$\hat{\mathbb G}(\Omega)$ is Banach \cite{Dav}, we have that the factor-algebra $\mathbb G(\Omega)$ is also Banach \cite{Gam}.
\end{proof}

\begin{lemma}
\label{teo2}
$\mathbb{PC}(\Omega)$ is dense in $\mathbb G(\Omega)$.
\end{lemma}

The proof is analogous to the proof of Lemma \ref{teo1}.

\begin{lemma}
\label{pointconv}
The map $\mathbb G(\Omega) \to \mathbb L^\infty(S^{n-1})$ given by $g \to g(x)(\cdot)$ is a continuous homomorphism.
\end{lemma}
\begin{proof}
As follows from the arithmetic properties of the limit, this map is a homomorphism. The continuity follows from (\ref{lem6conv}) (see proof of Lemma \ref{lem5}).
\end{proof}
Let us define the set of points of discontinuity of $g \in \mathbb G(\Omega)$ to be
\begin{equation*}
T(g)=\{x \in \Omega: g(x)(\cdot) \not\equiv \const \text{ in } \mathbb L^{\infty}(S^{n-1})\}. 
\end{equation*}
In what follows, by use of notation $g(x)$ we assume that $x \in \Omega \setminus T(g)$. 

\begin{lemma}
\label{Tlem}
Suppose that $g \in \mathbb G(\Omega)$. Then
$T(g) \subset \cup_{k=1}^\infty \partial S_k$
for certain $S_k \in \mathcal F$.
\end{lemma}
\begin{proof}
Let us denote by $\hat{T}(\hat{g})$ the set of points of discontinuity of a representative $\hat{g} \in \hat{\mathbb G}(\Omega)$. Clearly, $T(g) \subset \hat{T}(\hat{g})$. According to \cite{Dav} there exist
$S_k \in \mathcal F$ ($k \in \mathbb N$) such that $$\hat{T}(\hat{g})= \cup_{k=1}^\infty \partial S_k.$$ Consequently, $T(g) \subset \cup_{k=1}^\infty \partial S_k$.
\end{proof}
As it follows from Lemma \ref{Tlem}, for any $g \in \mathbb G(\Omega )$ the set of points of discontinuity $T(g)$ has zero measure, so the integration of $g \in \mathbb G(\Omega)$ is reduced to the integration of its representative in $\hat{\mathbb G}(\Omega)$. 
We also define the support of $g \in \mathbb G(\Omega)$ by the formula $$\supp(g)=\cl\{x \in \Omega: g(x)(\cdot)\ne 0 \text{ in } L^{\infty}(S^{n-1}) \}.$$

Let us denote by $\mathbb{PC}(\Omega)\hookrightarrow \mathbb G(\Omega)$ the image of the algebra of piecewise-constant functions under the canonical map $\hat{\mathbb G}(\Omega) \to \mathbb G(\Omega)$.
Also, let $\mathbb C(\Omega) \hookrightarrow \mathbb G(\Omega)$ be the algebra of continuous elements of $\mathbb G(\Omega)$.

\section{Distributions}

Let $\mathcal D(\Omega)$ be the space of test functions $\varphi \in \mathbb C(\Omega)$ having compact support $\supp(\varphi) \subset \Omega$ and endowed with the standard locally-convex topology (e.g., see \cite{Shi}).

Let $\mathcal R(\Omega)$ be the space of functions $\varphi \in \mathbb G(\Omega)$ having compact support $\supp(\varphi) \subset \rm \Omega$ with the fundamental system of neighbourhoods of zero $\{U_\gamma\}_{\gamma \in \mathbb C(\Omega),\gamma>0}$ in $\mathcal R'(\Omega)$ consisting of sets
\begin{equation}
U_\gamma=\{\varphi \in \mathcal R(\Omega)~:~|\varphi(x)(\cdot)|<\gamma(x)~(x \in \Omega)\},
\end{equation}
where
$|\varphi(x)(\cdot)|<\gamma(x)$ if $|\varphi(x)(s)|<\gamma(x)$ for all $s \in S^{n-1}$. 
As it follows straightforwardly from the definition of the topology in $\mathcal D(\Omega)$ \cite{Shi}, we have the embedding
\begin{equation*}
\mathcal D(\Omega) \hookrightarrow \mathcal R(\Omega).
\end{equation*}

\begin{theorem}
\label{lem_conv}
The space $\mathcal R(\Omega)$ is a locally-convex topological vector space.
\end{theorem}
\begin{proof}
Let us show that for any two neighbourhoods $U_{\gamma_1}$,
$U_{\gamma_2}$ there exists a neighbourhood
$U_{\gamma_3}$ such that
$U_{\gamma_3} \subset U_{\gamma_1} \cap U_{\gamma_2}$. 
Clearly, $U_{\gamma_1} \cap U_{\gamma_2}$ consists of the test functions $\varphi \in \mathcal R(\Omega)$ such that $|\varphi(x)(\cdot)|<\gamma_1(x)$, $|\varphi(x)(\cdot)|<\gamma_2(x)$ for any $x \in \Omega$. So, it suffices to put
\begin{equation*}
\gamma_3(x)=\min\{\gamma_1(x),\gamma_2(x)\} \quad x \in \Omega,
\end{equation*}
where $\gamma_3 \in \mathbb C(\Omega)$, $\gamma_3(x)>0$ for every  $x \in \Omega$ and $U_{\gamma_3}=U_{\gamma_1}\cap U_{\gamma_2}$.

Further, given $\lambda \in \mathbb R$, $|\lambda| \leqslant 1$, and a neighbourhood $U_\gamma$, we have
$\lambda U_\gamma \subset U_\gamma$
since for each $\varphi \in U_\gamma$
\begin{equation*}
|\lambda \varphi(x)(\cdot)|=|\lambda||\varphi(x)(\cdot)| \leqslant |\varphi(x)(\cdot)| <\gamma(x) 
\end{equation*}
for every $x \in \Omega$.

Suppose that $C \subset \Omega$ is compact.
Let $\varphi \in \mathcal R(\Omega)$, $\supp(\varphi) \subset C$, and $\gamma \in \mathbb C(\Omega)$, $\gamma>0$, be given.
Since $\min_{x \in C}\{|\gamma(x)|\}>0$, we may define
\begin{equation*}
\lambda=\frac{\max_{x \in C}\{\|\varphi(x)(\cdot)\|_{\mathbb L^{\infty}(S^{n-1})}\}}{\min_{x \in C}\{|\gamma(x)|\}} \geqslant 0.
\end{equation*}
Then, clearly,
$\varphi \in \mu U_{\gamma}$ 
for any $\mu \in \mathbb R$, $|\mu| \geqslant \lambda$.

Also, note that for any neighbourhood $U_\gamma$ there exists a neighbourhood $U_{\gamma'}$ such that $U_{\gamma'}+U_{\gamma'} \subset U_{\gamma}$. Indeed, we may put
$\gamma'=\gamma/2$. 

According to \cite{Kan}, the above results imply that
$\mathcal R(\Omega)$
is a topological vector space. 
Furthermore, consider the sequence of neighbourhoods $\{U_{\gamma_k}\}_{k=1}^\infty$, where
$\gamma_k \equiv 1/k$ on $\Omega$. 
Then 
$\cap_{k=1}^\infty U_{\gamma_k}=\{0\}$, so
according to \cite{Kan}
the space $\mathcal R(\Omega)$ is Hausdorff.

Finally, we note that given a neighbourhood $U_\gamma$ and the test functions
$\varphi, \psi \in U_\gamma$, we have 
\begin{equation*}
|\lambda \varphi(x)(\cdot)+(1-\lambda)\psi(x)(\cdot)| \leqslant
\lambda |\varphi(x)(\cdot)|+(1-\lambda)|\psi(x)(\cdot)|
\leqslant  \lambda \gamma(x)+(1-\lambda)\gamma(x)=\gamma(x), 
\end{equation*}
for all $x \in \Omega$ and $0 \leqslant \lambda \leqslant 1$, so $U_\gamma$ is convex. Since $\mathcal R(\Omega)$ is Hausdorff, this implies by definition that $\mathcal R(\Omega)$ is locally-convex.
\end{proof}

\begin{theorem}
\label{thmconv}
Let $\{\varphi_k\}_{k=1}^\infty \subset \mathcal R(\Omega)$, $\varphi \in \mathcal R(\Omega)$. Then $\varphi_k \to \varphi$ in $\mathcal R(\Omega)$ if and only if $\varphi_k \to \varphi$ in $\mathbb G(\Omega)$ and there exists a compact subset $C \subset \Omega$ such that $\supp(\varphi_k) \subset C$ \rm ($k \in \mathbb N$).
\end{theorem}

\begin{proof}
Suppose that there exists a compact subset $C \subset \Omega$ such that $\supp(\varphi_k) \subset C$ for all $k \in \mathbb N$. Also, without loss of generality we may assume that $\varphi_k \to 0$ in $\mathbb G(\Omega)$. As follows from the definition of the norm in $\mathbb G(\Omega)$, for any $\varepsilon>0$ there exists $N=N_\varepsilon \in \mathbb N$ such that 
\begin{equation}
\label{G_sup}
\|\varphi_k\|_{\mathbb G(\Omega)}<\varepsilon
\end{equation}
for all $k \geqslant N$. Let $U_\gamma$ be arbitrary. We define
$$\varepsilon=\inf_{x \in C}\{\gamma(x)/2\}>0.$$
Then $\varphi_k \in U_\gamma$ for all $k \geqslant N$ as follows from (\ref{G_sup}). Consequently, $\varphi_k \to 0$ in $\mathcal R(\Omega)$.

Now let $\varphi_k \to 0$ in $\mathcal R(\Omega)$. Let us consider the sequence of neighbourhoods $\{U_{\gamma_m}\}_{m=1}^\infty$, where $\gamma_m \equiv 1/m$. As follows from the definition of the convergence in $\mathcal R(\Omega)$, for any $m \in \mathbb N$ there exists $K \in \mathbb N$ such that $\varphi_k \in U_{\gamma_m}$ for any $k \geqslant K$, that is, $\varphi_k \to 0$
in $\mathbb G(\Omega)$. Suppose that there is no compact subset $C \subset \Omega$ such that $\supp(\varphi_k) \subset C$ for all $k \in \mathbb N$. Then there exist the sequences $\{x_k\}_{k=1}^\infty \subset \Omega$, $\{r_k\}_{k=1}^\infty \subset \mathbb R$ such that
\begin{equation*}
0<r_k\not<|\varphi(x_k)(\cdot)|, \quad r_k \to 0, 
\end{equation*}
and $x_k \to x$ in $\mathbb R^n$, where $x \in \partial \Omega$ or $\{x_k\}_{k=1}^\infty$ is unbounded. Without loss of generality we may consider the first case only and assume that $\{x_k\}_{k=1}^\infty$ does not have limit points in $\Omega$. Then there exists a function $\gamma \in \mathbb C(\Omega)$ such that $\gamma(x)>0$ ($x \in \Omega$) and $\gamma(x_k)=r_k$. Consequently, $\varphi_k \not\in U_\gamma$ for all $k \in \mathbb N$, i.e., $\varphi_k \not\to 0$ in $\mathcal R(\Omega)$. The contradiction obtained completes the proof. 
\end{proof}

%

Let $\mathcal D'(\Omega)$ and $\mathcal R'(\Omega)$ be the spaces of linear continuous functionals defined on $\mathcal D(\Omega)$ and $\mathcal R(\Omega)$, respectively.
The elements of both spaces $\mathcal D'(\Omega)$ and $\mathcal
R'(\Omega)$ are called the distributions.

\begin{theorem}
Any distributions in $\mathcal D'(\Omega)$ admits an extension from $\mathcal D(\Omega)$ to $\mathcal R(\Omega)$.
\end{theorem}
\begin{proof}
We have $\mathcal D(\Omega) \hookrightarrow \mathcal R(\Omega)$. According to Theorem \ref{lem_conv} the space $\mathcal R(\Omega)$ is locally-convex, so the extension exists by Hanh-Banach Theorem \cite{Kan}.
\end{proof}

The space $\mathcal R'(\Omega)$ is endowed with the linear operations and weak* topology, so by definition $f_k \to f$ in $\mathcal R'(\Omega)$ if and only if $(f_k,\varphi) \to (f,\varphi)$ for any $\varphi \in \mathcal R(\Omega)$. 
\begin{example}
\label{exreg}
Let $f \in \mathbb L_{\loc}(\Omega)$. Let us define the regular distribution $f \in \mathcal R'(\Omega)$ by
\begin{equation}
(f,\varphi)=\int_{\Omega} f(x)\varphi(x)dx, \quad \varphi \in \mathcal R(\Omega).
\end{equation}
Since $\mathcal D(\Omega) \hookrightarrow \mathcal R(\Omega)$ and the canonical map $\mathbb L_{\loc}(\Omega) \to \mathcal D'(\Omega)$ is injective \cite{Shi},
we may identify the elements of $\mathbb L_{\loc}(\Omega)$ and the regular distributions in $\mathcal R'(\Omega)$.
\end{example}

\begin{example}
\label{exdelta1}
Suppose that $p \in \Omega$. Suppose that we are given a function $$\alpha \in \mathbb L(S^{n-1}), \quad \int_{S^{n-1}} \alpha(s)ds=1.$$ We define the delta-function $\delta_p^\alpha$ by the formula
\begin{equation}
(\delta_p^\alpha,\varphi)=\int_{S^{n-1}}\alpha(s)\varphi(p)(s)ds. 
\end{equation}
The linearity and continuity of $\delta_\tau^\alpha$ follows from Lemma \ref{pointconv}, so $\delta_p^\alpha \in \mathcal R'(\Omega)$. 
Note that given any $\varphi \in \mathcal D(\Omega)$ we have that $\varphi(p)(\cdot) \equiv \varphi(p)$, so
\begin{equation*}
(\delta_p^\alpha,\varphi)=\int_{S^{n-1}}\alpha(s)\varphi(p)ds=\varphi(p).
\end{equation*}
Thus, $\delta_p^\alpha \in \mathcal R'(\Omega)$ is an extension of the classical delta-function $\delta_p \in \mathcal D'(\Omega)$ from $\mathcal D(\Omega)$ to $\mathcal R(\Omega)$.
\end{example}

\begin{example}
\label{onedeltaex}
Consider the case $n=1$. 
We define the right and the left delta-functions
$$(\delta_p^+,\varphi)=\varphi(p+), \quad (\delta_p^-,\varphi)=\varphi(p-), \quad \varphi \in \mathcal R(\Omega),$$
so
\begin{equation}
\label{oneform}
\delta_p^\alpha=\beta\delta_p^++(1-\beta)\delta_p^-,
\end{equation}
where $\beta=\alpha(1)$.
\end{example}

The support $\supp(f) \subset \Omega$ of $f \in \mathcal R'(\Omega)$ is the minimal closed set such that for any $\varphi \in \mathcal R(\Omega)$,
$\supp(f) \cap \supp(\varphi)=\varnothing$ we have $(f,\varphi)=0$. 
The distribution $f \in \mathcal R'(\Omega)$ is called non-negative {\rm (}we denote $f \geqslant 0${\rm)}, if $(f,\varphi) \geqslant 0$ for every $\varphi \geqslant 0$, $\varphi \in \mathcal R(\Omega)$. 

\begin{lemma}
\label{dist_complete}
If $\{f_k\}_{k=1}^\infty$ converges in $\mathcal R'(\Omega)$, $\varphi_k \to 0$ in $\mathcal R(\Omega)$,
then $(f_k,\varphi_k) \to 0$.
\end{lemma}

The proof of the following lemma is similar to the proof of analogous statement for the space $\mathcal D'(\Omega)$ in \cite{Shi}.

\begin{proof}
Suppose the contrary. Then without loss of generality (if necessary, let us consider a subsequence) we may suppose that there exists $c>0$ such that $|(f_k,\varphi_k)| \geqslant c>0$ for all $k \in \mathbb N$. Since $\varphi_k \to 0$ in $\mathcal R(\Omega)$, then we may suppose that
$\|\varphi_k\|_{\mathbb G(\Omega)} \leqslant \frac{1}{4^k}$. We define $$\zeta_k=2^k\varphi_k \in \mathcal R(\Omega).$$ Then we obtain the inequality
\begin{equation}
\label{oneineq100}
\|\zeta_k\|_{\mathbb G(\Omega)} \leqslant \frac{1}{2^k},
\end{equation}
so $\zeta_k \to 0$ in $\mathcal R(\Omega)$, but $$|(f_k,\zeta_k)|=2^k|(f_k,\varphi_k)| \geqslant 2^kc \to \infty.$$ We proceed further by induction. Let us choose $f_{k_1}$, $\zeta_{k_1}$ such that
$|(f_{k_1},\zeta_{k_1})|>1.$ Suppose that $f_{k_j}$, $\zeta_{k_j}$ were constructed ($1 \leqslant j \leqslant l-1$). Then for all $k \geqslant k'$ we have $$|(f_{k_j},\zeta_k)|<\frac{1}{2^{l-j}}, \quad 1 \leqslant j \leqslant l-1.$$
There exists $k_l \geqslant k'$ such that 
\begin{equation}
\label{oneineq101}
|(f_{k_l},\zeta_{k_l})|>\sum_{j=1}^{l-1}|(f_{k_l},\zeta_{k_j})|+l
\end{equation}
since $|(f_k,\zeta_k)| \to \infty$, $(f_k,\zeta_{k_j}) \to 0$ ($k \to \infty$). Suppose that the sequence $\{\zeta_{k_l}\}_{l=1}^\infty$ is constructed. Let us define
$\zeta=\sum_{j=1}^\infty \zeta_{k_j}$,
where the series converges due to (\ref{oneineq100}),
so $\zeta \in \mathcal R(\Omega)$. Consequently, 
\begin{equation*}
(f_{k_l},\zeta)=\sum_{j=1}^{l-1}(f_{k_l},\zeta_{k_j})+(f_{k_l},\zeta_{k_l})+\sum_{l+1}^\infty (f_{k_l},\zeta_{k_j}). 
\end{equation*}
Since (\ref{oneineq101}) and 
\begin{equation*}
\sum_{j=l+1}^\infty (f_{k_l},\zeta_{n_j})<\sum_{j=l+1}^\infty \frac{1}{2^{j-l}}=1, 
\end{equation*}
we obtain that $|(f_{k_l},\zeta)|>l-1$, which contradicts to the equality $(f_k,\zeta)=(f,\zeta)$, where $f=\lim_{k \to \infty} f_k$. The proof is complete.
\end{proof}


\begin{theorem}
The space of distributions $\mathcal R'(\Omega)$ is weak* complete.
\end{theorem}
\begin{proof}
We need to show that if $\{f_k\}_{k=1}^\infty \subset \mathcal R'(\Omega)$ and
$(f,\varphi)=\lim_{k \to \infty} (f_k,\varphi)$ for all $\varphi \in \mathcal R(\Omega)$, then $f \in \mathcal R'(\Omega)$.
Clearly, $f$ is a linear functional. In order to show the continuity of $f$, it suffices to prove that $f$ is continuous at $0 \in \mathcal R(\Omega)$.
So, let $\varphi_k \to 0$. Suppose that 
$(f,\varphi_k) \not\to 0$. Let $\varepsilon_0>0$ be given. Without loss of generality we may assume 
that for any $k \in \mathbb N$
\begin{equation*}
|(f,\varphi_k)|>\varepsilon_0.
\end{equation*}
As it follows from the definition, for any $k \in \mathbb N$ there exists
$m_k$ such that $|(f_{m_k},\varphi)|>\frac{\varepsilon_0}{2}$. Without loss of
generality we put $m_k=k$, so $|(f_k,\varphi_k)|>\frac{\varepsilon_0}{2}$ ($k \in \mathbb N$). This contradicts to the statement of Lemma \ref{dist_complete}. Thus, $(f,\varphi_k) \to 0$.
\end{proof}

Let us define the product of $f \in \mathcal R'(\Omega)$ and $g \in \mathbb G(\Omega)$ by the formula
\begin{equation}
(gf,\varphi)= (f,g\varphi),
\end{equation}
where $\varphi \in \mathcal R(\Omega)$, and $g\varphi \in \mathcal R(\Omega)$ since $\mathcal R(\Omega)$ is an ideal in $\mathbb G(\Omega)$. The operation of multiplication in $\mathcal R'(\Omega)$ is commutative and associative in the sense that for any $g$, $h \in \mathbb G(\Omega)$, $f \in \mathcal R'(\Omega)$ we have the identity
$$(gh)f=g(hf)$$
in $\mathcal R'(\Omega)$. 

\begin{example} 
\label{regmult}
Let $f \in \mathcal R'(\Omega)$ be a regular distribution, $g \in \mathbb G(\Omega)$. Then
\begin{equation*}
(gf,\varphi)= (f,g\varphi)=\int_\Omega f(x)g(x)\varphi(x)dx.
\end{equation*}
So, for regular distributions the operation of multiplication coincides with the ordinary one.
\end{example}

\begin{theorem}
\label{teomult}
Suppose that $g_k \to g$ in $\mathbb G(\Omega)$, $f_k \to f$ in $\mathcal R'(\Omega)$. Then $g_kf_k \to gf$ in $\mathcal R'(\Omega)$.
\end{theorem}
\begin{proof}
Let us note that $g_k \varphi \to g\varphi$ in $\mathcal R(\Omega)$
for any $\varphi \in \mathcal R(\Omega)$. Consequently, 
\begin{multline}
\notag
|(g_kf_k,\varphi)-(gf,\varphi)|=|(f_k,g_k\varphi)-(f,g\varphi)|
\leqslant \\ \leqslant |(f_k,g_k\varphi)-(f_k,g\varphi)|+|(f_k,g\varphi)-(f,g\varphi)|\leqslant\\
\leqslant |(f_k,g_k\varphi-g\varphi)|+|(f_k,g\varphi)-(f,g\varphi)| \to 0,
\end{multline}
in virtue of Lemma \ref{dist_complete} and due to convergence $f_k \to f$ in $\mathcal R'(\Omega)$.
\end{proof}

\begin{example}
Let $n=1$, $p \in \Omega$. The Heaviside function $\theta_p \in \mathbb G(\Omega)$ is defined by
\begin{equation*}
\theta_p(x)=\left\{
\begin{array}{ll}
1,&x>p,\\
0,&x<p.
\end{array}
\right.
\end{equation*} 
Then
$(\theta_p\delta_p^\alpha,\varphi)= (\delta_p^\alpha,\theta_p\varphi)=\beta\theta_p(p+)\varphi(p+)+(1-\beta)\theta_p(p-)\varphi(p-)=\beta\varphi(p+)$,
i.e.,
\begin{equation*}
\theta_p\delta_p^\alpha=\beta\delta_p^+,
\end{equation*}
where $\beta=\alpha(1)$ (see Example \ref{onedeltaex}).
\end{example}
\begin{example}
Let $n=2$, $p \in \Omega$. Let us find the product of the function $g \in \mathbb G(\Omega)$ given by
\begin{equation*}
g(x)=\left\{
\begin{array}{ll}
\sigma,& x^1>p^1,~x^2>p^2, \\
\mu, & x^1<p^1,~x^2<p^2, \\
0, & \text{ otherwise},
\end{array}
\right.
\end{equation*}
and the delta-function $\delta_p^\alpha \in \mathcal R'(\Omega)$, where $\alpha:[0,2\pi) \mapsto \mathbb R$, $\int_0^{2\pi}\alpha(s)ds=1$. We have
\begin{equation}
\label{gdelta}
(g\delta_p^\alpha,\varphi)=\int_{0}^{2\pi}\alpha(s)g(p)(s)\varphi(p)(s)ds=
\int_0^{\frac{\pi}{2}}\alpha(s)\varphi(p)(s)ds+\int_{\pi}^{\frac{3\pi}{2}}\alpha(s)\varphi(p)(s)ds
\end{equation}
Denote $\rho=\int_0^{\frac{\pi}{2}}\alpha(s)ds+\int_{\pi}^{\frac{3\pi}{2}}\alpha(s)ds \in \mathbb R$.
If $\rho \ne 0$, then the equality (\ref{gdelta}) can be rewritten as
\begin{equation}
\label{deltaprod1}
g\delta_p^\alpha=\rho\delta_p^\gamma,
\end{equation}
where $\gamma$ is  given by $\gamma(s)=\frac{\sigma \alpha(s)}{\rho}$ if $0<s<\frac{\pi}{2}$, $\gamma(s)=\frac{\mu \alpha(s)}{\rho}$ if $\pi<s<\frac{3\pi}{2}$,
$\gamma(s)=0$ otherwise.
\end{example}

Let us define the integral of $f \in \mathcal R'(\Omega)$ over $S \in \mathcal F_c = \{S \in \mathcal F: \cl(S) \subset \Omega\}$ by the formula
\begin{equation}
\label{int}
\int_S f dx = (f,\chi_S),
\end{equation}
where the characteristic function $\chi_S \in \mathcal R(\Omega)$. The integral (\ref{int}) exists for any distribution, is linear as a function of $f$ and coincides with the Lebesgue integral for regular distributions.

\begin{theorem}
\label{intconv}
Suppose that $S \in \mathcal F_c$, $f_k \to f$ in $\mathcal R'(\Omega)$. Then
\begin{equation*}
\int_S f_k dx \to \int_S f dx.
\end{equation*}
\end{theorem}
\begin{proof}
We have
$\int_S f_k dx=(f_k,\chi_S) \to (f,\chi_S)=\int_S f dx$.
\end{proof}

\begin{example}
Let $n=1$, $p$, $t_0 \in \Omega$, $t_0<p$. Then
\begin{equation*}
\int_{t_0}^t \delta_p^\alpha ds=\left\{
\begin{array}{ll}
1,&t>p, \\
1-\beta, & t=p, \\
0, &t<p.
\end{array}
\right.
\end{equation*}
\end{example}
\begin{example}
Let $n=2$, $p \in \Omega$, $B_p^2 \subset \Omega$ is a disk centered at $p$, $S_p(r) \subset B_p^2$ is a sector possessing central angle $r \in [0,2\pi)$. Then
\begin{equation}
\label{int2ex}
\int_{S_p(r)}\delta_p^\alpha dx= (\delta_p^\alpha,\chi_{S_p(r)})=\int_0^{2\pi}\alpha(s)\chi_{S_p(r)}(p)(s)ds=\int_0^r \alpha(s)ds, \quad r \in [0,2\pi),
\end{equation}
where $\chi_{S_p(r)}(p)(\cdot)$ is the surrounding value of the characteristic functions $\chi_{S_p(r)}$ at $p\in \Omega$. 
\end{example}

In what follows, we denote by $K_p \subset \mathbb R^n$ a closed convex cone with the vertex $p \in \Omega$ such that $K_p \cap \Omega \in \mathcal F$. Let $B^p_{r} \subset \Omega$ be an open ball of radius $r>0$ centered at $p$. Also, we put
\begin{equation*}
B_r(K_p)= K_p \cap B^p_r, \quad S(K_p) = (K_p-p) \cap S^{n-1}.
\end{equation*}
As it follows immediately from the definition, we have $B_r(K_p) \in \mathcal F$. 

The next theorem is the main result of the present section.

\begin{theorem}
\label{deltaseqteo}
Suppose that $\{\omega_k\}_{k=1}^\infty$ is a delta-sequence, that is, $\omega_k \to \delta_p^\alpha$ in $\mathcal R'(\Omega)$. Then
\begin{equation}
\label{conv1}
\int_{B_r(K_p)}\omega_k(x)dx \to \int_{S(K_p)}\alpha(s)ds.
\end{equation}
for any $K_p$ and any $r>0$.

Conversely, if $$\alpha \in \mathbb L^\infty(S^{n-1}), \quad \int_{S^{n-1}}\alpha(s)ds=1,$$ for any cone $K_p$ and any radius $r>0$ we have \rm(\ref{conv1}), for any $r>0$ we have
$w_k \to \delta_p^\alpha$ in $\mathcal R'(\Omega \setminus \cl(B_r^p))$,
and there exists $C>0$ such that for any two cones $K_p' \subset K_p''$
\begin{equation*}
\sup_{k \in \mathbb N}\int\limits_{B_r(K_p'') \setminus B_r(K_p')} |\omega_k(x)|dx \leqslant C {\rm mes }\bigl(S(K_p'') \setminus S(K_p') \bigr),
\end{equation*}
then 
$\omega_k \to \delta_p^\alpha$ in $\mathcal R'(\Omega)$.
\end{theorem}
\begin{proof}
1) There exists $r_0>0$ such that for any $0<r<r_0$ we have $\cl(B_r(K_p)) \subset \Omega$. Since $B_r(K_p) \in \mathcal F$, for any such $r>0$ we have $\chi_{B_r(K_p)} \in \mathcal R(\Omega)$ and 
\begin{equation*}
(\omega_k,\chi_{B_r(K_p)})=\int_{\Omega}\omega_k(x)\chi_{B_r(K_p)}(x)dx=\int_{B_r(K_p)}\omega_k(x)dx,
\end{equation*}
\begin{equation*}
(\delta_p^\alpha,\chi_{B_r(K_p)})=\int_{S^{n-1}}\alpha(s)\chi_{B_r(K_p)}(p)(s)=\int_{S(K_p)}\alpha(s)ds,
\end{equation*}
where the last equality follows from the fact that $K_p$ is a cone.
By definition, if $\omega_k \to \delta_p^\alpha$ in $\mathcal R'(\Omega)$ then
\begin{equation*}
\bigl(\omega_k,\chi_{B_r(C_p)}\bigr) \to \bigl(\delta_p^\alpha,\chi_{B_r(C_p)}\bigr), 
\end{equation*}
that is, we have (\ref{conv1}).

2) The convergence $\omega_k \to \delta_p^\alpha$ in $\mathcal R'(\Omega)$ is equivalent to
\begin{equation}
\label{reqconv1}
\int_\Omega \omega_k(x)\varphi(x)dx \to \int_{S^{n-1}}\alpha(s)\varphi(p)(s)ds
\end{equation}
for any $\varphi \in \mathcal R(\Omega)$. 
It suffices to show that (\ref{reqconv1}) is true for any $\varphi \in \mathcal R(\Omega) \cap \mathbb{PC}(\Omega)$. Indeed, for any $\varphi \in \mathcal R(\Omega)$ there exists a sequence $\{\varphi_l\}_{l=1}^\infty \subset \mathcal R(\Omega) \cap \mathbb{PC}(\Omega)$, $\varphi_l \to \varphi$ in $\mathcal R(\Omega)$, i.e., $$\varphi_l \to \varphi \text{ in } \mathbb G(\Omega), \quad \varphi_l(p)(\cdot) \to \varphi(p)(\cdot) \text{ in } \mathbb L^\infty(S^{n-1})$$ according to Lemma \ref{pointconv}, and
\begin{equation*}
\label{reqconv2}
\lim_{l \to \infty}\lim_{k \to \infty} \int_\Omega \omega_k(x)\varphi_l(x)dx \to \lim_{l \to \infty} \int_{S^{n-1}}\alpha(s)\varphi_l(p)(s)ds=\int_{S^{n-1}}\alpha(s)\varphi(s)ds.
\end{equation*}
Further, given $l$, $m \in \mathbb N$, $k \in \mathbb N$, we have
\begin{equation*}
\biggl|\int_\Omega \omega_k(x)\varphi_l(x)dx-\int_\Omega \omega_k(x)\varphi_m(x)dx \biggr| 
\leqslant \int_\Omega |\omega_k(x)||\varphi_l(x)-\varphi_m(x)|dx \leqslant 4\pi C \|\varphi_l-\varphi_m\|_{\mathbb G(\Omega)},
\end{equation*}
so $$\int_\Omega \omega_k(x)\varphi_l(x)dx \to \int_\Omega \omega_k(x)\varphi(x)dx$$ uniformly with respect to $k \in \mathbb N$. Consequently, we may change the order of limits in (\ref{reqconv2}) to obtain (\ref{reqconv1}).

So, we need to show that (\ref{reqconv1}) holds for any $\varphi \in \mathcal R(\Omega) \cap \mathbb{PC}(\Omega)$. Without loss of generality we may restrict our consideration to the case $\varphi=\chi_M$ ($\cl(M) \subset \Omega$, $M$ is convex, i.e., $M \in \mathcal F$). In the special case if $M=B_r(K_p)$, where $K_p$ is a closed convex cone, then
\begin{equation}
\label{reqconv25}
\int\limits_\Omega \omega_k(x)\chi_{B_r(K_p)}(x)dx=\int\limits_{B_r(K_p)}\omega_k(x)dx \to \int\limits_{S(K_p)}\alpha(s)ds=\int\limits_{S^{n-1}}\alpha(s)\chi_{B_r(K_p)}(p)(s)ds
\end{equation}
according to the assumption of the theorem, i.e., we have (\ref{reqconv1}). In the general case for $M$ above we define a tangent cone $K_p(M)$ at $p$, $K_p(M)=\cl(C_p(M))$, where (see \cite{Aub2})
\begin{equation}
\label{conerepr}
C_p(M)=\{p+th:t>0,~h \in \mathbb R^n, \text{ where } p+ht \in M,~0<t<t_0(h)\}. 
\end{equation}
Since $K_p(M)$ is closed and convex \cite{Aub2}, we have (\ref{reqconv25}) for $K_p=K_p(M)$. Show that
\begin{equation}
\label{reqconv3}
\int_{S(K_p(M))}\alpha(s)ds=\int_S \alpha_S \chi_M(p)(s)ds,
\end{equation}
\begin{equation}
\label{reqconv4}
\int_{B_r(K_p(M))}\omega_k(x)dx-\int_\Omega \omega_k(x)\chi_M(x)dx \to 0,
\end{equation}
then the convergence (\ref{reqconv1}) for $\varphi=\chi_M$ will follow from (\ref{reqconv25}), (\ref{reqconv3}) and (\ref{reqconv4}).

Let us show that (\ref{reqconv3}) is true. Suppose that $s \in S^{n-1}$. Due to convexity of $M$ we have that $\chi_M(p)(s)=1$ if and only if there exists $t_0 \in (0,1)$ such that for any $0<t<t_0$ the point $p+ts \in M$ \cite{Aub2}. Consequently, there exists $S' \subset S^{n-1}$ such that
\begin{equation*}
\chi_M(p)(\cdot)=\chi_{S'}(\cdot):S^{n-1} \to \mathbb R. 
\end{equation*}
Let $C_p$ be the cone corresponding to $p$ and $S'$. It follows from (\ref{conerepr}) that $K_p(M)=\cl(C_p)$ and
\begin{equation*}
\int_{S(K_p(M))}\alpha(s)ds=\int_{S(C_p)}\alpha(s)ds=\int_{S'}\alpha(s)ds=\int_{S^{n-1}}\alpha(s)\chi_M(p)(s)ds,
\end{equation*}
i.e., the equality (\ref{reqconv3}) is true.

Let us show that (\ref{reqconv4}) is valid. We have that
\begin{equation*}
\int_\Omega \omega_k(x)\chi_M(x)dx=\int_{B_r^p}\omega_k(x)\chi_M(x)dx+\int_{\Omega \setminus B_r^p}\omega_k(x)\chi_M(x)dx,
\end{equation*}
where the last summand tends to 0 since $\omega_k \to 0$ in $\mathcal R'(\Omega \setminus \cl(B_{r'}^p))$ if $r'<r$. Then for a given $\varepsilon>0$ there exist $r_0>0$ and $K \in \mathbb N$ such that for any $0<r<r_0$ and any $k \geqslant K$
\begin{equation*}
\left|\int_{\Omega \setminus B_r^p}\omega_k(x)\chi_M(x)dx \right|<\varepsilon.
\end{equation*}
According to the remark above, we have to show that
\begin{equation}
\label{reqconv5}
\int_{B_r(K_p(M))}\omega_k(x)dx-\int_{B_r^p}\omega_k(x)\chi_M(x)dx \to 0.
\end{equation}
Observe that
\begin{equation*}
\biggl|\int_{B_r(K_p(M))}\omega_k(x)dx-\int_{B_r^p}\omega_k(x)\chi_M(x)dx \biggr| \leqslant 
\int_{B_r^p}|\omega_k(x)||\chi_{B_r(K_p(M))}(x)-\chi_M(x)|dx.
\end{equation*}
Further, $M \cap B_r^p \subset B_r(K_p)$ \cite{Aub2}. There is a cone $K_p' \subset K_p(M)$ such that
\begin{equation*}
\mes\bigl(S\bigl(K_p(M)\bigr) \setminus S(K_p')\bigr)<\varepsilon
\end{equation*}
and $B_r(K_p') \subset M \cap B_r^p$ for all $r>0$ sufficiently small. Consequently,
\begin{equation*}
|\chi_{B_r(K_p(M))}(x)-\chi_M(x)| \leqslant |\chi_{B_r(K_p(M))}(x)-\chi_{B_r(K_p')}(x)|,
\end{equation*}
where $x \in \Omega$.
Then the conditions of our theorem imply that
\begin{multline}
\notag
\biggl|\int_{B_r(K_p(M))}\omega_k(x)dx-\int_{B_r^p}\omega_k(x)\chi_M(x)dx \biggr| \leqslant 
\int_{B_r^p}|\omega_k(x)||\chi_{B_r(K_p(M))}-\chi_{B_r(K_p')}(x)|dx= \\=
\int_{B_r(K_p(M)) \setminus B_r(K_p')}|\omega_k(x)| 
\leqslant C  \mes\bigl(S(K_p(M)) \setminus S(K_p')\bigr)<C\varepsilon
\end{multline}
for all $k \in \mathbb N$ sufficiently large and all $r>0$ sufficiently small.
Consequently, we have convergence (\ref{reqconv5}) and, as a result, convergence (\ref{reqconv4}). 
\end{proof}

\begin{example}
\label{deltaseqex}
Let $n=1$, $\beta \in \mathbb R$. Then according to Theorem \ref{deltaseqteo} the sequence $\{\omega_k\}_{k=1}^\infty$,
\begin{equation*}
\omega_k = k(\beta \chi_{\left(p,p+\frac{1}{k}\right)}+(1-\beta)\chi_{\left(p-\frac{1}{k},p\right)})
\end{equation*}
is a delta-sequence, that is,
$\omega_k \to \delta_p^\alpha$,
where $\beta=\alpha(1)$.
\end{example}

\section{Zero-sum games with discontinuous payoff functions}

Let $\Omega=\Omega_1 \times \Omega_2 \subset \mathbb R^2$, where $\Omega_1$, $\Omega_2$ are open intervals in $\mathbb R$. Consider the following zero-sum game:
\begin{equation}
\label{game}
G=(X_1,X_2,\rho),
\end{equation}
where $X_1$, $X_2$ are the open intervals in $\mathbb R$ such that $\cl(X_1) \subset \Omega_1$, $\cl(X_2) \subset \Omega_2$. 
\subsection{Pure strategies and classical mixed strategies} The elements $x_1 \in X_1$ and $x_2 \in X_2$ such that $(x_1,x_2) \not\in T(\rho)$ are called the pure strategies of the first and the second player, respectively, the function
$\rho \in \mathbb G(\Omega)$ is called the payoff function of the first player \cite{Kan}.

Let us also consider $G$ in the set of mixed strategies, that is, the game $G^L=(X_1^L,X_2^L,\rho^L)$,
\begin{equation*}
X_1^L=\left\{u_1 \in \mathbb L(X_1): u_1 \geqslant 0, \int_{X_1}u_1(x_1)dx_1=1\right\}, 
\end{equation*}
\begin{equation*}
X_2^L=\left\{u_2 \in \mathbb L(X_2): u_2 \geqslant 0, \int_{X_2}u_2(x_2)dx_2=1\right\} 
\end{equation*}
-- the sets of mixed strategies of the first and the second player, respectively, the map 
\begin{equation*}
\rho^L(u_1,u_2) = \int_{X_1 \times X_2} \rho(x_1,x_2)u_1(x_1)u_2(x_2)dx_1dx_2. 
\end{equation*}
is called the payoff function of the first player. 
As the following example shows, $G$ may have no solution in the sets of pure or mixed strategies.

\begin{example}
\label{exgame}
Let $X_1=X_2=(-1,1)$, we define
\begin{equation*}
\rho(x_1,x_2)=\left\{
\begin{array}{ll}
1,&x_1,x_2>0,x_1+x_2<1 \text{ or } \\
&x_1,x_2<0,x_1+x_2>-1, \\
0, & \text{otherwise}.
\end{array}
\right.
\end{equation*}
Let us consider game 
$G=(X_1,X_2,\rho)$. 

First, observe that $G$ does not have a solution in the set of pure strategies $x_1 \in X_1 \setminus \{0\}$,
$x_2 \in X_2 \setminus \{0\}$. Indeed, for any $x_2$ we have $\sup_{x_1}\rho(x_1,x_2)=1$,
i.e., 
$\inf_{x_2}\sup_{x_1}\rho(x_1,x_2)=1$. 
Similarly, for any $x_1$ we have $\inf_{x_2}\rho(x_1,x_2)=0$, so
$\sup_{x_1}\inf_{x_2}\rho(x_1,x_2)=0$. 
According to \cite{Kan} $G$ does not have a solution.

Second, let us show that $G$ does not have solution in the set of mixed strategies. We define
\begin{equation*}
\sigma_{u_1}(x_2)=\int_{X_1}\rho(x_1,x_2)u_1(x_1)dx_1,
\end{equation*}
where $u_1 \in X_1^L$. Then $\sigma_{u_1} \geqslant 0$, function $\sigma_{u_1}$ is monotonically increasing on $(-1,0) \subset X_2$ and is monotonically decreasing on $(0,1) \subset X_2$, $\sigma_{u_1}(1-)=\sigma_{u_1}(-1)=0$. 
For any $\varepsilon>0$ there exists $u^\varepsilon_2 \in X_2^L$, $\supp(u^\varepsilon_2) \subset (1-\varepsilon,1)$, such that
\begin{equation*}
\rho^L(u_1,u_2^\varepsilon)=\int_{X_2}\sigma_{u_1}(x_2)u^\varepsilon_2(x_2)dx_2<\varepsilon.
\end{equation*}
Consequently, $\inf_{u_2}\rho^L(u_1,u_2)=0$ for any $u_1 \in X_1^L$. Then
\begin{equation}
\label{e0}
\sup_{u_1}\inf_{u_2}\rho^L(u_1,u_2)=0.
\end{equation}
Analogously, for a given $u_2 \in X_2^L$ we define
\begin{equation*}
\tau_{u_2}(x_1)=\int_{X_2}\rho(x_1,x_2)u_2(x_2)dx_2.
\end{equation*}
Then $\tau_{u_2} \geqslant 0$ is monotonically increasing on $(-1,0) \subset X_1$,  is decreasing on $(0,1) \subset X_1$, and $\tau_{u_2}(0+)+\tau_{u_2}(0-)=1$.
For any $\varepsilon>0$ there exists $u_1^\varepsilon \in X_1^L$, $\supp(u_1^\varepsilon) \subset (-\varepsilon,\varepsilon)$, such that
\begin{equation*}
\rho^L(u_1^\varepsilon,u_2)=\int_{X_1}\tau_{u_2}(x_1)u_1^\varepsilon(x_1)dx_1>1-\varepsilon.
\end{equation*}
Consequently, $\sup_{u_1}\rho^L(u_1,u_2)=1$ for any $u_2 \in X_2^L$. Then
\begin{equation}
\label{e1}
\inf_{u_2}\sup_{u_1}\rho^L(u_1,u_2)=1.
\end{equation}
Now comparison of (\ref{e0}) and (\ref{e1}) shows that $G^L$ does not have solution (see \cite{Kan}).
\end{example}

\subsection{$\mathcal R'$-mixed strategies}In order to provide the existence of solution, let us consider game $G$ in the set of $\mathcal R'$-mixed strategies, that is, the game $G^R=(X_1^R,X_2^R,\rho^R)$, where the sets $X_1^R$ and $X_2^R$ consist of distributions
$v_1 \in \mathcal R'(\Omega_1)$ and $v_2 \in \mathcal R'(\Omega_2)$ such that
\begin{equation*}
v_1 \geqslant 0, \quad \int_{X_1}v_1dx_1=1, \quad v_2 \geqslant 0, \quad \int_{X_2}v_2dx_2=1,
\end{equation*}
the functions
\begin{equation*}
x_1 \to \int_{X_2}\rho(x_1,x_2)v_2dx_2 \in \mathbb G(\Omega_1), \quad x_2 \to \int_{X_1}\rho(x_1,x_2)v_1dx_1
\end{equation*}
are in $\mathbb G(\Omega_2)$
and the following equality is satisfied:
\begin{equation}
\label{game2}
\int_{X_1}\! \biggl(\int_{X_2}\!\rho(x_1,x_2)v_2dx_2\biggr)v_1dx_1\!=\!\int_{X_2}\! \biggl(\int_{X_1}\!\rho(x_1,x_2)v_1dx_1\biggr)v_2dx_2=:\rho^R(v_1,v_2)
\end{equation}
(the definitions of a non-negative distribution and the integral of a distribution were given earlier).

The sets $X_1^R$ and $X_2^R$ are called the sets of $\mathcal R'$-mixed strategies of the first and the second player, respectively, the map $\rho^R$ is called the payoff function of the first player.

Let us note that $X_1^L \subset X_1^R$, $X_2^L \subset X_2^R$ (see Example \ref{exreg}), $\rho^R|_{X_1^L \times X_2^L}=\rho^L$, so the classical mixed strategies for $G$ can be viewed as $\mathcal R'$-mixed strategies. 
We denote
\begin{equation*}
\rho(x_1^* \pm,x_2) = \lim_{x_1 \to x_1^* \pm} \rho(x_1,x_2), \quad \rho(x_1 ,x_2^*\pm) = \lim_{x_2 \to x_2^* \pm} \rho(x_1,x_2).
\end{equation*}

\begin{theorem}
\label{gameteo}
Let $\rho \geqslant 0$.
If there exists $(x_1^*,x_2^*) \in X_1 \times X_2$ such that 
\begin{equation}
\label{a1}
a_{\pm}^r = \lim_{x_1 \to x_1^*+}\rho(x_1,x_2^*\pm) \geqslant \rho(x_1,x_2^* \pm), \quad x_1>x_1^*,
\end{equation}
\begin{equation}
\label{a2}
a_{\pm}^l = \lim_{x_1 \to x_1^*-}\rho(x_1,x_2^*\pm) \geqslant \rho(x_1,x_2^* \pm), \quad x_1<x_1^*,
\end{equation}
\begin{equation}
\label{b1}
b_{\pm}^r = \lim_{x_2 \to x_2^*+}\rho(x_1^*\pm,x_2) \leqslant \rho(x_1^*\pm,x_2), \quad x_2>x_2^*,
\end{equation}
\begin{equation}
\label{b2}
b_{\pm}^l = \lim_{x_2 \to x_2^*-}\rho(x_1^*\pm,x_2) \leqslant \rho(x_1^*\pm,x_2), \quad x_2<x_2^*,
\end{equation}
\begin{equation}
\label{equa}
b_+^r=a_+^r, \quad b_-^l=a_-^l, \quad b_-^r=a_+^l, \quad b_+^l=a_-^r,
\end{equation}
\begin{equation}
\label{ni}
a_+^r \geqslant a_-^r, \quad a_-^l \geqslant a_+^l, \quad a_-^l \geqslant a_-^r, \quad a_+^r \geqslant a_+^l, 
\end{equation}
and $a_+^r-a_+^l \ne a_-^r-a_-^l$, then the pair of delta-functions 
$\delta_{x_1^*}^{\alpha^*_1} \in X_1^R$, $\delta_{x_2^*}^{\alpha^*_2} \in X_2^R$, 
where
\begin{equation*}
\alpha^*_1(1)=\frac{a_-^l-a_+^l}{a_+^r-a_+^l-a_-^r+a_-^l}, \quad \alpha^*_2(1)=\frac{a_-^l-a_-^r}{a_+^r-a_+^l-a_-^r+a_-^l},
\end{equation*}
gives the solution of game $G$ in the set of $\mathcal R'$-mixed strategies.
\end{theorem}
Since the function $\rho$ is bounded on $\Omega$, the case $\rho \not\geqslant 0$ can be reduced to the one considered above by consideration of the payoff function $\rho+C$ for $C>0$ sufficiently large.
\begin{proof}
1) Let us show that $\delta_{x_1^*}^{\alpha^*_1}$, $\delta_{x_2^*}^{\alpha^*_2}$ are $\mathcal R'$-mixed strategies, that is, the value $\rho(\delta_{x_1^*}^{\alpha^*_1},\delta_{x_2^*}^{\alpha^*_2})$ is correctly defined.
Denote
$A=a_+^r-a_+^l-a_-^r+a_-^l \ne 0$.
Observe that
\begin{equation*}
\int_{X_2}\rho(x_1,x_2)\delta_{x_2^*}^{\alpha^*_2}dx_2=\frac{a_-^l-a_-^r}{A}\rho(x_1,x_2^*+)+\biggl(1-\frac{a_-^l-a_-^r}{A}\biggr)\rho(x_1,x_2^*-) \in \mathbb G(\Omega_1).
\end{equation*}
Consequently, we have the following equality:
\begin{multline}
\notag
\int_{X_1}\!\biggl(\int_{X_2}\rho(x_1,x_2)\delta_{x_2^*}^{\alpha^*_2}dx_2\biggr)\delta_{x_1^*}^{\alpha^*_1}dx_1\!=\!\frac{a_-^l-a_+^l}{A}\biggl(\frac{a_-^l-a_-^r}{A}a_+^r+\biggl(1-\frac{a_-^l-a_-^r}{A}\biggr)a_-^r\!\biggr)+\\
+\biggl(1-\frac{a_-^l-a_+^l}{A}\biggr)\biggl(\frac{a_-^l-a_-^r}{A}a_+^l+\biggl(1-\frac{a_-^l-a_-^r}{A}\biggr)a_-^l \biggr)=\frac{a_+^ra_-^l-a_-^ra_+^l}{A}.
\end{multline}
Analogously,
\begin{equation*}
\int_{X_1}\rho(x_1,x_2)\delta_{x_1^*}^{\alpha^*_1}dx_1=\frac{a_-^l-a_+^l}{A}\rho(x_1^*+,x_2)+\biggl(1-\frac{a_-^l-a_+^l}{A}\biggr)\rho(x_1^*-,x_2) \in \mathbb G(\Omega_2).
\end{equation*}
According to (\ref{equa}), we have
\begin{multline}
\notag
\int_{X_2}\!\biggl(\int_{X_1}\rho(x_1,x_2)\delta_{x_1^*}^{\alpha^*_1}dx_1\biggr)\delta_{x_2^*}^{\alpha^*_2}dx_2\!=\!\frac{a_-^l-a_-^r}{A}\biggl(\frac{a_-^l-a_+^l}{A}b_+^r+\biggl(1-\frac{a_-^l-a_+^l}{A}\biggr)b_-^r\!\biggr)+\\
+\biggl(1-\frac{a_-^l-a_-^r}{A}\biggr)\biggl(\frac{a_-^l-a_+^l}{A}b_+^l+\biggl(1-\frac{a_-^l-a_+^l}{A}\biggr)b_-^l \biggr)=\frac{a_+^ra_-^l-a_-^ra_+^l}{A}.
\end{multline}
Then by our definition $\delta_{x_1^*}^{\alpha^*_1}$, $\delta_{x_2^*}^{\alpha^*_2}$ are $\mathcal R'$-mixed strategies for $G$. 

2) Observe that we always have the following inequality:
\begin{equation}
\label{kani}
\inf_{v_2}\sup_{v_1}\rho^R(v_1,v_2) \geqslant \sup_{v_1}\inf_{v_2}\rho^R(v_1,v_2)
\end{equation}
(see \cite{Kan}). Let $X_2=(p_1,p_2) \subset \mathbb R$. Then we obtain
\begin{multline}
\label{oneeq}
\rho(\delta_{x_1^*}^{\alpha^*_1},v_2)=\int_{X_2}\biggl(\int_{X_1}\rho(x_1,x_2)\delta_{x_1^*}^{\alpha^*_1}dx_1\biggr)v_2dx_2=\\
=
\int_{X_2}\biggl(\frac{a_-^l-a_+^l}{A}\rho(x_1^*+,x_2)+\biggl(1-\frac{a_-^l-a_+^l}{A}\biggr)\rho(x_1^*-,x_2) \biggr) v_2dx_2
\end{multline}
Further, according to the definitions of the product and the integral, 
\begin{multline}
\label{ras}
\int_{p_1}^{x_2^*}\rho(x_1^*+,x_2)v_2dx_2=\int_{p_1}^{x_2^*}b_+^lv_2dx_2+
\int_{p_1}^{x_2^*}(\rho(x_1^*+,x_2)-b_+^l)v_2dx_2=\\=\int_{p_1}^{x_2^*} b_+^l v_2dx_2+(v_2,(\rho(x_1^*+,x_2)-b_+^l)\chi_{(p_1,x_2^*)}) \geqslant \int_{p_1}^{x_2^*} b_+^l v_2dx_2
\end{multline}
since $v_2 \geqslant 0$ in $\mathcal R'(\Omega_2)$, $\rho(x_1^*,x_2) \geqslant b_+^l$ ($x<x_2^*$) according to (\ref{b2}). Since $\rho \geqslant 0$, $A \ne 0$ and (\ref{ni}) holds, we obtain that $A>0$. Then (\ref{b1})--(\ref{equa}), the equality (\ref{oneeq}) and the argument similar to (\ref{ras}) gives us the inequality
\begin{equation*}
\begin{split}
\rho(\delta_{x_1^*}^{\alpha^*_1},v_2)
\geqslant 
\int_{p_1}^{x_2^*}&\biggl(\frac{a_-^l-a_+^l}{A}b_+^l+\biggl(1-\frac{a_-^l-a_+^l}{A}\biggr)b_-^l\biggr)v_2dx_2+\\+
\int_{x_1^*}^{p_2}&\biggl(\frac{a_-^l-a_+^l}{A}b_+^r+\biggl(1-\frac{a_-^l-a_+^l}{A}\biggr)b_-^r\biggr)v_2dx_2 = \\
=
\int_{p_1}^{x_2^*}&\biggl(\frac{a_-^l-a_+^l}{A}a_-^r+\biggl(1-\frac{a_-^l-a_+^l}{A}\biggr)a_-^l\biggr)v_2 dx_2+\\
+\int_{x_2^*}^{p_2}&\biggl(\frac{a_-^l-a_+^l}{A}a_+^r+\biggl(1-\frac{a_-^l-a_+^l}{A}\biggr)a_+^l\biggr)v_2 dx_2
=\frac{a_+^ra_-^l-a_+^la_-^r}{A}.
\end{split}
\end{equation*}
for any $\mathcal R'$-mixed strategy $v_2$. Consequently,
\begin{equation*}
\inf_{v_2}\rho^R(\delta_{x_1^*}^{\alpha^*_1},v_2) \geqslant \frac{a_+^ra_-^l-a_+^la_-^r}{A}, \text{ i.e., } \sup_{v_1}\inf_{v_2}\rho^R(v_1,v_2) \geqslant \frac{a_+^ra_-^l-a_+^la_-^r}{A}.
\end{equation*}
Analogously, due to the inequalities (\ref{a1}), (\ref{a2}) and the equalities (\ref{equa}), we have 
\begin{equation*}
\inf_{v_2}\sup_{v_1}\rho^R(v_1,v_2) \leqslant \frac{a_+^ra_-^l-a_+^la_-^r}{A}. 
\end{equation*}
Then (\ref{kani}) implies that 
\begin{equation}
\label{game3}
\max_{v_1}\inf_{v_2}\rho^R(v_1,v_2)=\min_{v_2}\sup_{v_1}\rho^R(v_1,v_2), 
\end{equation}
i.e., the solution of $G^R$ exists, the maximum and the minimum 
in (\ref{game3}) 
are attained at
\begin{equation}
\label{game4}
v_1^*=\delta_{x_1^*}^{\alpha^*_1} \text{ and } v_2^*=\delta_{x_2^*}^{\alpha^*_2}, 
\end{equation}
respectively. According to \cite{Kan} the pair (\ref{game4}) forms a solution of $G^R$.
\end{proof}

\begin{example}
Let us consider game $G=(X_1,X_2,\rho)$ of Example \ref{exgame}.
Let $(x_1^*,x_2^*)=(0,0)$. Then the conditions of Theorem \ref{gameteo} are satisfied, where
\begin{equation*}
a_+^r=1, \quad a_+^l=0, \quad a_-^r=0, \quad a_-^l=1, \quad
b_+^r=1, \quad b_+^l=0, \quad b_-^r=0, \quad b_-^l=1,
\end{equation*}
so the pair
\begin{equation}
\label{sol}
\delta_{0}^{\alpha^*_1} \in X_1^R, \quad \delta_{0}^{\alpha^*_2} \in X_2^R, \quad \alpha^*_1(1)=\alpha^*_2(1)=\frac{1}{2}
\end{equation}
is a solution of $G$ in the $\mathcal R'$-mixed strategies.

Let us note that solution (\ref{sol}) admits an approximation by the mixed strategies (see Example \ref{deltaseqex}) and, thus, possesses an obvious probabilistic interpretation. 
\end{example}

\bibliographystyle{alpha}
\bibliography{derkin}

\end{document}